\documentstyle[11pt,mssymb]{article}

\font \twlmath = msym10 at 11pt
\newfam \mathfam \textfont \mathfam = \twlmath

\newcommand{\arrow}{\longrightarrow}

\newcommand{\Z}{{\Bbb Z}}
\newcommand{\C}{{\Bbb C}}
\newcommand{\R}{{\Bbb R}}

\newcommand{\g}{{\goth g}}

\newcommand{\1}{\sqrt{-1}\:}

\def\blacksquare{\hbox{\vrule width 4pt height 4pt depth 0pt}}

\begin{document}

\centerline{\bf Hyperk\"ahler embeddings}
\centerline{\bf and holomorphic symplectic geometry \  \rm I.}

\centerline{Mikhail Verbitsky,}
\centerline{verbit@math.harvard.edu}

\hfill

\hfill

{\bf 0. Introduction.}

 In this paper we are studying
complex analytic subvarieties of a given K\"ahler manifold
which is endowed with a holomorphic symplectic structure.

By Calabi-Yau theorem, the holomorphically symplectic K\"ahler manifolds
can be supplied with a  Ricci-flat Riemannian metric.
This implies that such manifolds
are hyperk\"ahler (Definition 1.1). Conversely, all hyperk\"ahler
manifolds are holomorphically symplectic (Proposition 2.1).

For a given closed analytic subvariety $S$ of a holomorphically symplectic
$M$, one can restrict the holomorphic symplectic form of $M$
to the Zarisky tangent sheaf to $S$. If this restriction
is non-degenerate outside of singularities of $S$,
this subvariety is called
non-degenerately symplectic. (Definition 2.2).
Of course, such subvarieties are of even complex dimension.

Take a generic element $N$ in a given deformation class of a
holomorphically symplectic K\"ahler manifolds.
We are proving that all complex analytic subvarieties of $N$ are
non-degenerately symplectic (Theorem 2.3). In particular,
all closed analytic subvarieties of $N$ are of even complex dimension.
If such subvariety is smooth, it is also a hyperkaehler manifold
(Proposition 2.1).

\hfill

{\sf \centerline{Contents.}}

{\sf 1. Hyperk\"ahler manifolds.}

{\sf 2. Holomorphic symplectic geometry.}

{\sf 3. The action of $\goth{so}(5)$ on the differential forms over a
hyperk\"ahler manifold.}

\hfill

In the first section we give basic definitions from the theory
of the hyperk\"ahler manifolds ([Bes]). In the second section,
we state our results and deduce them from Theorem 2.1.
In the third section, we give the proof of Theorem 2.1.
This section is based entirely on calculations from [V].

\hfill

{\bf 1. Hyperk\"ahler manifolds.}

\hfill

{\bf Definition 1.1} ([B], [Bes]) A {\bf hyperk\"ahler manifold} is a
Riemannian manifold $M$ endowed with three complex structures $I$, $J$
and $K$, such that the following holds.

\hspace{5mm}   (i)  $M$ is K\"ahler with respect to these structures and

\hspace{5mm}   (ii) $I$, $J$ and $K$, considered as  endomorphisms
of a real tangent bundle, satisfy the relation
$I\circ J=-J\circ I = K$.

\hfill

This means that the hyperk\"ahler manifold has the natural action of
quaternions ${\Bbb H}$ in its real tangent bundle.
Therefore its complex dimension is even.


Let $\mbox{ad}I$, $\mbox{ad}J$ and $\mbox{ad}K$ be the operators on the
bundles of differential forms over a hyperk\"ahler manifold
$M$ which are defined as follows. Define $\mbox{ad}I$.
Let this operator act as a complex structure operator
$I$ on the bundle of differential 1-forms. We
extend it on $i$-forms for arbitrary $i$ using Leibnitz
formula: $\mbox{ad}I(\alpha\wedge\beta)=\mbox{ad}I(\alpha)\wedge\beta+
\alpha\wedge \mbox{ad}I(\beta)$. Since Leibnitz
formula is true for a commutator in a Lie algebras, one can immediately
obtain the following identities, which follow from the same
identities in ${\Bbb H}$:

\[
   [\mbox{ad}I,\mbox{ad}J]=2\mbox{ad}K;\;
   [\mbox{ad}J,\mbox{ad}K]=2\mbox{ad}I;\;
\]

\[
   [\mbox{ad}K,\mbox{ad}I]=2\mbox{ad}J
\]

Therefore, the operators $\mbox{ad}I,\mbox{ad}J,\mbox{ad}K$
generate a Lie algebra $\goth{su}(2)$ acting on the
bundle of differential forms. We can integrate this
Lie algebra action to the action of a Lie group
$SU(2)$. In particular, operators $I$, $J$
and $K$, which act on differential forms by the formula
$I(\alpha\wedge\beta)=I(\alpha)\wedge I(\beta)$,
belong to this group.

{\bf Proposition 1.1:} There is an action of the Lie group $SU(2)$
and Lie algebra $\goth{su}(2)$ on the bundle of differential
forms over a hyperk\"ahler manifold. This action is
parallel, and therefore it commutes with Laplace operator.
$\blacksquare$

If $M$ is compact, this implies that there is
a canonical $SU(2)$-action on $H^i(M,\R)$ (see [V1]).

\hfill

Let $M$ be a hyperk\"ahler manifold with a Riemannian form $<\cdot,\cdot>$.
Let the form $\omega_I := <I(\cdot),\cdot>$ be the usual K\"ahler
form  which is closed and parallel
(with respect to the connection). Analogously defined
forms $\omega_J$ and $\omega_K$ are
also closed and parallel.

The simple linear algebraic
consideration ([B]) shows that \hfill
$\omega_J+\sqrt{-1}\omega_K$ is of
type $(2,0)$ and, being closed, this form is also holomorphic.
It is called {\bf the canonical holomorphic symplectic form
of a manifold M}. Conversely, if there is a parallel
holomorphic symplectic form on a K\"ahler manifold $M$,
this manifold has a hyperk\"ahler structure ([B]).

If some $compact$ K\"ahler manifold $M$ admits non-degenerate
holomorphic symplectic form $\Omega$, the Calabi-Yau ([Y]) theorem
implies that $M$ is hyperk\"ahler (Proposition 2.1).
This follows from the existence of a K\"ahler
metric on $M$ such that $\Omega$ is parallel for the Levi-Civitta
connection associated with this metric.

\hfill

Let $M$ be a hyperk\"ahler manifold with complex structures
$I$, $J$ and $K$. For any real numbers $a$, $b$, $c$
such that $a^2+b^2+c^2=1$ the operator $L:=aI+bJ+cK$ is also
an almost complex structure: $L^2=-1$.
Clearly, $L$ is parallel with respect to a connection.
This implies that $L$ is a complex structure, and
that $M$ is K\"ahler with respect to $L$.

\hfill

{\bf Definition 1.2} If $M$ is a  hyperk\"ahler manifold,
the complex structure $L$ is called {\bf induced
by a hyperk\"ahler structure}, if $L=aI+bJ+cK$ for some
real numbers $a,b,c\:|\:a^2+b^2+c^2=1$.

\hfill

\hfill

If $M$ is a hyperk\"ahler manifold and $L$ is induced complex structure,
we will denote $M$, considered as a complex manifold with respect to
$L$, by $(M,L)$ or, sometimes, by $M_L$.

\hfill

Consider the Lie algebra $\goth{g}_M$ generated by ${ad}L$ for all $L$
induced by a hyperk\"ahler structure on $M$. One can easily see
that $\goth{g}_M=\goth{su}(2)$.
The Lie algebra $\goth{g}_M$ is called {\bf isotropy algebra} of $M$, and
corresponding Lie group $G_M$ is called an {\bf isotropy group}
of $M$. By Proposition 1.1, the action of the group is parallel,
and therefore it commutes with Laplace operator in differential
forms. In particular, this implies that the action of the isotropy
group $G_M$ preserves harmonic forms, and therefore this
group canonically acts on cohomology of $M$.

\hfill

{\bf Proposition 1.2:} Let $\omega$ be a differential form over
a hyperk\"ahler manifold $M$. The form $\omega$ is $G_M$-invariant
if and only if it is of Hodge type $(p,p)$ with respect to all
induced complex structures on $M$.

{\bf Proof:} Assume that $\omega$ is $G_M$-invariant.
This implies that all elements of $\g_M$ act trivially on
$\omega$ and, in particular, that $\mbox{ad}L(\omega)=0$
for any induced complex structure $L$. On the other hand,
$\mbox{ad}L(\omega)=(p-q)\1$ if $\omega$ is of Hodge type $(p,q)$.
Therefore $\omega$ is of Hodge type $(p,p)$ with respect to any
induced complex structure $L$.

Conversely, assume that $\omega$ is of type $(p,p)$ with respect
to all induced $L$. Then $\mbox{ad}L(\omega)=0$ for any induced $L$.
By definition, $\g_M$ is generated by such $\mbox{ad}L(\omega)=0$,
and therefore $\g_M$ and $G_M$ act trivially on $\omega$. $\blacksquare$

\hfill

\hfill

{\bf 2. Holomorphic symplectic geometry.}

\hfill

{\bf Definition 2.1:} The compact K\"ahler manifold $M$ is called
holomorphically symplectic if there is a holomorphic 2-form $\Omega$
over $M$ such that $\Omega^n=\Omega\wedge\Omega\wedge...$ is
a nowhere degenerate section of a canonical class of $M$.
There, $2n=dim_\C(M)$.

Note that we assumed compactness of $M$.
\footnote{If one wants to define a holomorphic symplectic
structure in a situation when $M$ is not compact,
one should require also the equation $\nabla'\Omega$ to held.
The operator $\nabla':\;\Lambda^{p,0}(M)\arrow\Lambda^{p+1,0}(M)$
is a holomorphic differential defined on differential forms ([GH]).}
One sees that the holomorphically symplectic
manifold has a trivial canonical bundle.
A hyperk\"ahler manifold is holomorphically symplectic
(see Section 1). There is a converse proposition:

{\bf Proposition 2.1} ([B], [Bes])  Let $M$ be a holomorphically
symplectic K\"ahler manifold with the holomorphic symplectic form
$\Omega$, a K\"ahler class
$[\omega]\in H^{1,1}(M)$ and a complex structure $I$.
There is a unique hyperk\"ahler structure $(I,J,K,(\cdot,\cdot))$
over $M$ such that the cohomology class of the symplectic form
$\omega_I=(\cdot,I\cdot)$ is equal to $[\omega]$ and the
canonical symplectic form $\omega_J+\1\omega_K$ is
equal to $\Omega$.

Proposition 2.1 immediately
follows from the Calabi-Yau theorem ([Y]). $\:\;\blacksquare$

\hfill

For each complex analytic variety $X$ and a point $x\in X$,
we denote the Zariski tangent space to $X$ in $x$ by $T_xX$.

{\bf Definition 2.2:} Let $M$ be a holomorphically symplectic
manifold and $S\subset M$ be its complex analytic subvariety.
Assume that $S$ is closed in $M$ and reduced. It is called
non-degenerately symplectic if for each point $s\in S$
outside of the singularities of $S$ the restriction
of the holomorphic symplectic form $\Omega$ to $T_sM$
is nondegenerate on $T_s S\subset T_s M$, and the set
$Sing(S)$ of the singular points of $S$ is nondegenerately
symplectic. This definition refers to itself, but
since $dim\:Sing(S)<dim\:S$, it is consistent.

Of course, the complex dimension of a non-degenerately symplectic
variety is even.

\hfill





\hfill

Let $M$ be a holomorphically symplectic K\"ahler manifold.
By Proposition 2.1, $M$ has a unique hyperk\"ahler metric with
the same K\"ahler class and holomorphic symplectic form.
Therefore one can without ambiguity speak
about the action of $G_M$ on $H^*(M,\R)$ (see Proposition 1.1).
Of course, this action essentially depends on the choice
of K\"ahler class.

\hfill

{\bf Definition 2.3:} The K\"ahler form over a holomorphically symplectic
manifold $M$ is called {\bf of general type} if all elements of
$H^{pp}(M)\cap H^{2p}(M,\Z)$ are $G_M$-invariant, where the
action of $G_M$ is defined by Proposition 2.1. The holomorphically
symplectic manifold $M$ is called {\bf of general type} if there
exists a K\"ahler form of general type over $M$.

\hfill

As Theorem 2.2 implies, the holomorhically symplectic manifold
of general type has no Weil divisors. Therefore these manifolds
have connected Picard group. In particular, such manifolds
are never algebraic.

\hfill

{\bf Proposition 2.2:} Let $M$ be a hyperk\"ahler manifold. Let $S$
be the set of induced complex structures over $M$. Let $S_0\subset S$
be the set of $R\in S$ such that the natural K\"ahler metric
on $(M,R)$ is of general type. Then $S_0$ is dense in $S$.

{\bf Proof:}
Let $A$ be the set of all $\alpha\in H^{2p}(M,\Z)$ such that $\alpha$
is not $G_M$-invariant. The set $A$ is countable. For each
$\alpha\in A$, let $S_\alpha$ be the set of all $R\in S$
such that $\alpha$ is of type $(p,p)$ with respect to $R$.
The set $S_0$ of all induced complex structures of general type
is equal to $\{S\backslash\bigcup_{\alpha\in A}S_\alpha\}$.
Now, to prove Proposition 2.2 it is enough to prove that
$S_\alpha$ is a finite set for each $\alpha\in A$.
This would imply that $S_0$ is a complement of a
countable set to a 2-sphere $S$, and therefore dense in $S$.

As it follows from Section 1,
$\alpha$ is of type (1,1) with respect to $R$ if and
only if $ad \:R(\alpha)=0$. Now, let $V$ be a representation
of $\goth{su}(2)$, and $v\in V$ be a non-invariant vector.
It is easy to see that the element $a\in \goth{su}(2)$
such that $a(v)=0$ is unique up to a constant, if it exists.
This implies that if $\alpha$ is not $G_M$-invariant
there are no more than two $R\in S$ such that
$ad R(\alpha)=0$.
Of course, these two elements
of $S$ are opposite to each other. $\;\blacksquare$

\hfill

One can easily deduce from the results in [Tod] and
from Proposition  2.2 that the set of points
associated with holomorphically symplectic
manifolds of general type is dense in the classifying space
of holomorphically symplectic manifolds.



\hfill

\hfill

For a K\"ahler manifold $M$ and a form
$\alpha\in H^{2p}(M,\C)$, define
\[deg(\alpha):=\int_M L^p(\alpha)\]
where $L$ is a Hodge operator of exterior
multiplication by the K\"ahler form $\omega$
(see [GH]).
Of course, the degree of forms of Hodge type
$(p,q)$ with $p\neq q$ is equal zero,
so only $(p,p)$-form can possibly have
non-zero degree.

\hfill

We recall that the real dimension of a
holomorphically symplectic manifold is divisible by 4.

{\bf Theorem 2.1:} Let $M$ be a holomorhically
symplectic K\"ahler manifold with
a holomorphic symplectic form $\Omega$.
Let $\alpha$ be a $G_M$-invariant form of non-zero degree.
Then the dimension of $\alpha$ is divisible by 4.
Moreover,
\[ \int_M \Omega^n\wedge\bar\Omega^n\wedge\alpha=2^n deg(\alpha),\]
where $n=\frac{1}{4}(dim_R M-dim\:\alpha)$.

\hfill

This theorem will be proven in Section 3.

\hfill

{\bf Theorem 2.2:} Let $M$ be a holomorphic symplectic
manifold of general type. All closed analytic subvarieties
of $M$ have even complex dimension.

{\bf Proof:} Take a closed reduced analytic subvariety $S\subset M$.
Assume it represents a cycle $[S]$ in cohomology of $M$,
and take a cocycle $\alpha\in H^{2p}(M,\Z)$ which is dual to $[S]$
by Poincare duality. We have to prove that $p$ is even.

By definition of degree,
\[ deg(\alpha)=\int_S\omega\wedge\omega\wedge...\]
Since $S$ is complex analytic,
$\int_S\omega\wedge\omega\wedge...>0$ and therefore
$deg(\alpha)>0$.
On the other hand, since the form $\alpha$ is integer, it is
$G_M$-invariant because $M$ is of general type.
Therefore, by Theorem 2.1, $dim(\alpha)$ is divisible by 4, and
$dim_\C S=dim_\C M-\frac{1}{2}dim\: \alpha$ is even.
$\:\:\:\blacksquare$

\hfill

{\bf Theorem 2.3:} Let $M$ be a holomorphic symplectic
manifold of general type. All reduced closed analytic
subvarieties of $M$ are non-degenerately symplectic.

{\bf Proof:} Let $S\in M$ be a closed analytic subvariety of
$M$. Take the restriction
of $\Omega$ to $S$.  The complex dimension
of $S$ is even by Theorem 2.2. Let $dim_\C S=2p$
and $dim_\C M=2n$. By definition,

\[ \int_S(\Omega\wedge\bar\Omega)^p=
   \int_M(\Omega\wedge\bar\Omega)^p\wedge\alpha
\]
where $\alpha$ is a Poincare dual form to $[S]$ (see the proof
of Theorem 2.2).

By Theorem 2.1, the last integral is non-zero.
Therefore $\int_S(\Omega\wedge\bar\Omega)^p\neq 0$
and the holomorphic $p$-form $\Omega^p$ is a
nontrivial section of a canonical line bundle over $S$.
Take $Z$ to be the set of zeros of $\Omega^p$.
Clearly, $Z$ is a closed complex analytic
subvariety of $S\subset M$.
Outside of singularities of $S$,
the subvariety $Z\subset S$ has codimension 1 in $S$.
Since by Theorem 2.2 the holomorphic manifold
$M$ has no odd-dimensional closed analytic subsets,
$Z$ belongs to the set of singularities of $S$. Therefore the restriction
of $\Omega$ to $S$ is non-degenerate outside of
singularities. Applying the same consideration to
 $Sing(S)$, we see that $S$ is nondegenerately symplectic.
$\;\;\blacksquare$

\hfill

Combining this with Proposition 2.1, one obtains

{\bf Theorem 2.4:} Let $M$ be a holomorphically symplectic
manifold of general type, and $S\subset M$ be its smooth differentiable
submanifold. If $S$ is analytic in $M$, it is a
hyperkaehler manifold. $\;\;\blacksquare$

\hfill

\hfill

{\bf 3. The action of $\goth{so}(5)$ on the differential forms over a
hyperk\"ahler manifold.}

\hfill

In this section, we denote the space of smooth differential $i$-forms
over a manifold $M$ by $A^i(M)$. The notation for the
Hodge decomposition is $A^i(M)=\oplus_{p+q=i}A_R^{p,q}(M)$,
where $R$ is a complex structure operator the decomposition
is defined in respect with.

Let $M$ be a hyperk\"ahler manifold. For every induced complex structure
$R$ over $M$, there is a real symplectic form $\omega_R=(\cdot,R\cdot)$
(see Section 1). As usually, $L_R$ denotes the operator of exterior
multiplication by $\omega_R$, which is acting on the differential
forms $A^*(M,\C)$ over $M$. The operator of interior multiplication
by $\omega_R$, which is defined as an adjoint operator to $L_R$,
is denoted by $\Lambda_R$.

One may ask oneself, what algebra is generated by $L_R$ and
$\Lambda_R$ for all induced $R$? The answer was partially given
in [V], where the following theorem was proven.

\hfill

{\bf Theorem 3.1} ([V]) Let $M$ be a hyperk\"ahler manifold
and $\goth{a}_M$ be a Lie algebra generated by $L_R$ and $\Lambda_R$
for all induced complex structures $R$ over $M$. The
Lie algebra  $\goth{a}_M$ is isomorphic to $\goth{so}(5)$.
$\;\;\blacksquare$

\hfill

The following facts about a structure of $\goth{a}_M$
were proven in [V].

Let $I$, $J$ and $K$ be three induced complex structures on
$M$, such that $I\circ J=-J\circ I=K$
The algebra  $\goth{a}_M$ is 10-dimensional. It contains
$\g_M$ as a subalgebra, as follows:

\[ [\Lambda_J,L_K]=[L_J,\Lambda_K]= ad\: I\;\; {\rm (etc)}.\]
It has a following base: $L_R,\Lambda_R$, $ad\: R$
$(R=I,J,K)$ and the element $H=[L_R\Lambda_R]$.
Of course, $H$ is a standard Hodge operator,
which is independent on $R$. In acts on $r$-forms
over $M$ as a  multiplication by a scalar $n-r$,
where $n=dim_\C M$.

\hfill

The semisimple Lie algebra $\goth{a}_M$ has a two-dimensional
Cartan subalgebra $\goth{h}_M$, spanned over $\C$ by
$H$ and $ad\: I$. This algebra has a root system
$B_2$; the elements $H$ and $\1 ad\: I$ are
among its roots. This implies that the weight
decomposition of the $\goth{a}_M$-module
$A^*(M)$ taken with respect to this particular Cartan
subalgebra $\goth{h}_M\subset\goth{a}_M$ coincides
with the Hodge decomposition
$A^*(M)=\oplus_{p,q}A_I^{p,q}(M)$

\hfill

\hfill

{\bf Proof of Theorem 2.1$\;$} Let $\alpha$ be a $\g_M$-invariant
$2p$-form over a hyperkaehler manifold $M$
such that $\int_M L_I^{n-p}(\alpha)\neq 0$,
where $n=dim_\C(M)$. We have to prove that $p$ is even,
and that
\[ \int_M \Omega^{(n-p)/2}\wedge\bar\Omega^{(n-p)/2}\wedge\alpha
=2^{(n-p)/2}\int_M L_I^{n-p}(\alpha).\]

\hfill

Let $A^*(M)=\oplus_{l\in\Pi}A_l$ be the isotypic
decomposition of a $\goth a_M$-module
$A^*(M)$. Isotypic decomposition is defined as follows.
For each $l\in \Pi$, where $\Pi$ is a weight lattice of $\goth{a}_M$,
the module $A_l$ is a union of all simple $\goth{a}_M$-submodules
of $A^*(M)$ with a highest weight $l$.
One can easily see that the isotypic decomposition does not depend
on a choice of a Cartan subalgebra of $\goth{a}_M$. This follows,
for example, from Schuhr's lemma.

\hfill

Take the Cartan subalgebra $\goth{h}_M\subset\goth{a}_M$.
One can choose the simple roots of $\goth{h}_M$
to be $H$ and $\1 ad\; I-H$ in $\goth{h}_M$.
The vectors $H$, $H+\1 ad\; I$, $\1 ad\; I$
and $\1 ad\; I-H$ are the positive roots,
and $\rho= \1 ad\; I+H$ is the longest root.
The weight decomposition of $A^*(M)$
with respect to this particular Cartan algebra
is a Hodge decomposition.

{\bf Lemma 3.1:} If a simple $\goth{a}_M$-submodule $F\subset A^*(M)$
contains a $2n$-form $\theta$, then $\theta$
is a highest vector of $F$. We use the notation
$dim_\C M= n$.

{\bf Proof:} Choose Cartan subalgebra and the
positive roots as above. A $(p,q)$-form in $A^{p,q}_I(M)$
has a weight $p+q$ with respect to the longest root $\rho$.
Therefore the weight of the $(n,n)$-form is higher
that the weight of any $(p,q)$-form with $(p,q)\neq(n,n)$.
$\;\;\blacksquare$

\hfill

One of the summands of the isotypic decomposition is of
particular interest to us. We denote
this summand by $A_o$. The $\goth{a}_M$-module
$A_o\subset A^*(M)$ is a union of all simple
$\goth{a}_M$-submodules of  $A^*(M)$ which contain
a $2n$-form. By Lemma 3.1, $A_o$ is a summand of
the isotypic decomposition.

Consider the decomposition
$\displaystyle \alpha=\sum_{l\in\Pi}\alpha_l$,
which corresponds to the isotypic decomposition
of $A^*(M)$. All components $\alpha_l$
of this decomposition are $\g_M$-invariant, because
$\g_M\subset \goth{a}_M$.
The form $L_I^{n-p}(\alpha_l)$ has a dimension $2n$
for all $l\in\Pi$.
By definition, $L_I^{n-p}(\alpha_l)\in A_l$. Lemma 3.1
states that if $A_l$ contains a $2n$-form, there is a
simple submodule of $A_l$ which
has a highest weight $o$. Therefore if $L_I^{n-p}(\alpha_l)\neq 0$
then $l=o$.

\hfill

This implies that
$L_I^{n-p}(\alpha)=L_I^{n-p}(\alpha_o)$.

Assume that $p$ is even. By definition,
\[
  \Omega^{(n-p)/2}\wedge\bar\Omega^{(n-p)/2}\wedge\alpha=
  (L_J+\1L_K)^{(n-p)/2}(L_J-\1L_K)^{(n-p)/2}(\alpha)
\]

By the same reason as before,
\[
  (L_J+\1L_K)^{(n-p)/2}(L_J-\1L_K)^{(n-p)/2}(\alpha) =
\]
\[
  (L_J+\1L_K)^{(n-p)/2}(L_J-\1L_K)^{(n-p)/2}(\alpha_o),
\]
because $(L_J+\1L_K)^{(n-p)/2}(L_J-\1L_K)^{(n-p)/2}(\alpha)$
is of dimension $2n$. This implies that
\[
  \Omega^{(n-p)/2}\wedge\bar\Omega^{(n-p)/2}\wedge\alpha=
  \Omega^{(n-p)/2}\wedge\bar\Omega^{(n-p)/2}\wedge\alpha_o
\]

Therefore one can assume that $\alpha\in A_o$.
The following proposition proves Theorem 3.1.

\hfill

Take the restriction of $\alpha$ to the
fiber $\Lambda^*(T_xM)$ of $A^*(M)$ in the point $x\in M$.
Clearly, $\goth a_M$ and $\g_M$ act on $\Lambda^*(T_xM)$,
and the restriction of $\alpha$ to $T_xM$ is $\g_M$-invariant.

{\bf Proposition 3.1:} Let $A$ be a $\goth{a}_M$-submodule
of $\Lambda^*(T_xM)$ generated over $\goth{a}_M$ by
the unique determinant form $det\in\Lambda^{2n}(T_xM)$
and $\alpha\in A$ be a non-zero
$\g_M$-invariant form of dimension $2p$.
Then $p$ is even, and

\[\Omega^{(n-p)/2}\wedge\bar\Omega^{(n-p)/2}\wedge\alpha=
2^{(n-p)/2}L^{n-p}\alpha\]

{\bf Proof:} Let $S$ be the set of all induced complex structures over
$M$. Clearly, $S$ is a 2-sphere (see Section 1). If we identify
$G_M$ with $SU(2)$, then $S=\{s\in G_M\:|\: s^2=-1\}$,
where $-1$ is the matrix $-Id\in SU(2)$. For each
$R_1$, $R_2\in \{I, J, K\}$

\[ [ad R_1, \Lambda_{R_2}]=0 \;\;{\rm if}\; R_1=R_2\]
and
\[
   [ad R_1, \Lambda_{R_2}]=L_{R_1\circ R_2}
   \;\;\rm if\;\it R_1\neq R_2
\]
(see [V]).

Therefore, since $ad \:R (det)=0$ for each $R\in S$,
any form $\alpha\in A$ can be represented as a polynomial

\[ \alpha = P(\Lambda_I,\Lambda_J,\Lambda_K) det \]
with coefficients in $\C$. This representation is not unique:
there is an epimorphism map

\[ \C[\Lambda_I,\Lambda_J,\Lambda_K]\;\arrow\hspace{-0.7cm}\arrow A.\]

Denote the kernel of this map by $\goth R$. We have
$\C[\Lambda_I,\Lambda_J,\Lambda_K]/\goth R=A$.

The group $G_M$ acts on $S$ by the adjoint action.
For $t\in G_M$ and an operator $B$ on $\Lambda^*(T_xM)$,
denote $t\circ B\circ t^{-1}$ by $B^t$.
One can easily prove that $L_R^t=L_{R^t}$.
Since $\alpha$ is $G_M$-invariant, the
polynomial $P$ is $G_M$-invariant in the following sense:
\[ P(\Lambda_I^t,\Lambda_J^t,\Lambda_K^t)det =
P(\Lambda_I,\Lambda_J,\Lambda_K) det.\]

Let $V$ be a 3-dimensional space spanned by $\Lambda_J$, $\Lambda_J$
and $\Lambda_K$. Since $\forall R\in S$ the operator
$\Lambda_R$ is uniquely represented as a linear combination of
$\Lambda_J$, $\Lambda_J$ and $\Lambda_K$, we can consider $\Lambda_R$
as an element of $V$.
The group $G_M=SU(2)$ acts on $V$ as follows:
an element $t\in G_M$ maps $\Lambda_R$ in $\Lambda_{R^t}$.
Actually, his acion is a standard  action of
$SO(3)=SU(2)/\{\pm 1\}$ on a 3-dimensional
vector space.

The ring $\C[V]$ is graded by the degree of polynomials.
Clearly, $G_M$ preserves $\goth R$. Moreover, $\goth R$
is a graded ideal in a graded ring $\C[V]$, and the isomorphism
$A=\C[V]/\goth R$ is an isomorphism of graded $G_M$-modules.

Let $S^kV\subset \C[V]$ be the space of homogeneous polynomials
of degree $k$ over $V$ and $\goth R_k=S^kV\cap \goth R$.
Of course, $\oplus_{k}\goth R_k=\goth R$. Since $G_M$ is reductive,
one can represent $S^kV=\goth R_k\oplus \goth P_n$, where
$\goth P_n$ is a $G_M$-invariant submodule of $\C[V]$.
Therefore, $A$ can be embedded in $\C[V]$ as a $G_M$-module:

\[A:\:\stackrel{i}{\hookrightarrow} \C[V] \]
such that the composition

\[A:\:\stackrel{i}{\hookrightarrow} \C[V]\arrow \C[V]/\goth R=A\]
is identity.

Take the element $P=i(\alpha)$; it is a polynomial
of three variables $\Lambda_I$, $\Lambda_J$ and $\Lambda_K$
such that $\alpha = P(\Lambda_I^I,\Lambda_J^I,\Lambda_K^I) det$.
Moreover, $P$ is a $G_M$-invariant element of $\C[V]$,
because $\alpha$ is $G_M$-invariant.

Now, one can apply a simple theorem of the representation theory.

{\bf Lemma 3.2.} ([W])  If $SO(3)$ acts on a three-dimensional
space with an orthogonal basis $x,\:y,\:z$, then
all $SO(3)$-invariant polynomials in
$\C[V]$ are proportional to $C^k$ for some integer $k$,
where $C=x^2+y^2+z^2$. $\;\;\blacksquare$

\hfill

This lemma implies that
\[\alpha= (\Lambda_I^2+\Lambda_J^2+\Lambda_K^2)^k(det)\]
Therefore, $dim(\alpha)$ is divisible by 4. Moreover,
a simple calculation proves that

\[
   (L_J+\1L_K)^{(n-p)/2}(L_J-\1L_K)^{(n-p)/2}
   (\Lambda_I^2+\Lambda_J^2+\Lambda_K^2)^{(n-p)/2}(det)=
\]
\[
   = (L_J^2+L_K^2)^{(n-p)/2}
   (\Lambda_I^2+\Lambda_J^2+\Lambda_K^2)^{(n-p)/2}(det) =
\]
\[
   = 2^{(n-p)/2}L_I^{n-p}
   (\Lambda_I^2+\Lambda_J^2+\Lambda_K^2)^{(n-p)/2}(det)
\]

Theorem 3.1 is proven. $\;\;\blacksquare$

\hfill

I am very grateful to my advisor David Kazhdan for a warm support
and encouragement. Thanks also due to Michael Finkelberg
and Roman Bezrukavnikov for interesting discussions.

\hfill

\centerline{\bf Reference:}

\hfill

\hfill

[B] Beauville, A. Varietes K\"ahleriennes dont la pere classe de Chern est
nulle. // J. Diff. Geom. 18, p. 755-782 (1983).

\hfill

[Bes] Besse, A., Einstein Manifolds. // Springer-Verlag, New York (1987)

\hfill

[GH] Griffits, Ph. and Harris, J. Principles of algebraic geometry. //
Wiley-Interscience, New York (1978).

\hfill

[Tod] Todorov, A. Moduli of Hyper-K\"ahlerian manifolds I,II. // Preprint
MPI (1990)

\hfill

[V] Verbitsky, M. On the action of a Lie algebra SO(5) on the cohomology
of a hyperk\"ahler manifold. // Func. Analysis and Appl. 24(2)
p. 70-71 (1990).

\hfill

[W] Weyl, H. The classical groups, their invariants and representations.
// Princeton Univ. Press, New York (1939)

\hfill

[Y] Yau, S. T. On the Ricci curvature of a compact K\"ahler manifold
and the complex Monge-Amp\`ere equation I. // Comm. on Pure and Appl.
Math. 31, 339-411 (1978).

\end{document}